\documentclass[11pt]{amsart}   
\usepackage{amssymb,xy}
\usepackage{amsthm, amsmath}
\usepackage[latin1]{inputenc}
\usepackage[T1]{fontenc}
\usepackage{fancyvrb}
\usepackage{color}
\usepackage{url}
\xyoption{all}
\theoremstyle{plain}

\newtheorem{theorem}{Theorem}[section]

\theoremstyle{definition}

\newtheorem{definition}[theorem]{Definition}

\newtheoremstyle{algorithm}{3ex}{\topsep}{}{}{\bfseries}{.}{\newline}{}
\theoremstyle{algorithm}
\newtheorem{Alg}{Algorithm}

\theoremstyle{remark}

\newcommand{\inp}{{\itshape Input: }}
\newcommand{\outp}{{\itshape Output: }}
\newcommand{\lrrdc}{\mbox{\,\(\Leftarrow\hspace{-9pt}\Leftarrow\hspace{-5pt}\Rightarrow\hspace{-9pt}\Rightarrow\)\,}} 

\newcommand{\bmp}{\rule{0pt}{1pt}\\[-15pt]{\tiny.\dotfill.}\\*[-3pt]}
\newcommand{\bmpi}{\begingroup\footnotesize}
\newcommand{\empi}{\endgroup\rule{0pt}{1pt}\\*[-3pt]}
\newcommand{\empim}{\endgroup\hspace{5pt}\mbox{\small\(\maltese\)}\rule{0pt}{1pt}\\*[-3pt]}

\newcommand{\emp}{\rule{0pt}{1pt}\\*[-17pt]{\tiny.\dotfill.}\\*[-10pt]\rule{0pt}{1pt}}
\newcommand{\boxtt}[1]{\mbox{\small\texttt{#1}}}

\title{Effective Homology of the Pushout of Simplicial Sets}
\author{Jónathan Heras}
\date{}

\address{\small \rm  Departamento de Matemáticas y Computación, Universidad de La Rioja,
Edificio Vives, Luis de Ulloa s/n, E-26004 Logroño (La Rioja, Spain).}
\email{jonathan.heras@unirioja.es}

\begin{document}

\begin{abstract}
In this paper, an algorithm building the effective homology version of the pushout of the simplicial morphisms $f:X\rightarrow Y$ and $g:X\rightarrow Z$, where $X,Y$ and $Z$ are simplicial sets with effective homology is presented.
\end{abstract}

\maketitle

\section*{Introduction}

Many of the usual constructions in Topology are nothing but homotopy pullbacks or homotopy pushouts \cite{Mat76}. Loop spaces, suspensions, mapping cones, wedges or joins, for instance, can involve such constructions. In these cases, when the spaces are not of finite type the computation of their homology groups can be considered as a challenging task.

Sergeraert's ideas about effective homology \cite{RS02} provide real algorithms for the computation of homology groups of complicated spaces, such as loop spaces, classifying spaces, total spaces of fibrations and so on. The effective homology method has been concretely implemented in the Kenzo system \cite{Kenzo}, a successful Common Lisp program devoted to Symbolic Computation in Algebraic Topology which has obtained some homology groups which have not been confirmed nor refuted by any other means.

In this work we use the effective homology method to design algorithms building the pushout associated with two simplicial morphisms $f:X\rightarrow Y$ and $g:X\rightarrow Z$, where $X,Y$ and $Z$ are simplicial sets with \emph{effective homology}.

\section{Effective Homology of the Pushout}

\begin{definition}
A homotopy commutative diagram
$$\xymatrix{X\ar[r]^{f}\ar[d]^{g}&Y\ar[d]^{f'}&\\
Z\ar[r]^{g'}&P&\\
}$$
\noindent equipped with $H: f'\circ f\sim g'\circ g$, is called a \emph{homotopy pushout} \cite{Mat76}, denoted $P_{(f,g)}$, when for any commutative diagram
$$\xymatrix{X\ar[r]^{f}\ar[d]^{g}&Y\ar[d]^{f''}&\\
Z\ar[r]^{g''}&Q&\\
}$$
\noindent equipped with $G: f''\circ f\sim g''\circ g$, the following properties hold:
\begin{enumerate}
\item there exists a map $p:P\rightarrow Q$ and homotopies $K:f''\sim p\circ f'$ and $L: p\circ g' \sim g''$ such that the whole diagram $$\xymatrix{X\ar[r]^{f}\ar[d]^{g}&Y\ar[d]^{f'}\ar[rrd]^{f''}&\\
Z\ar[r]^{g'}\ar@{->}@/_/[rrr]_{g''}&P\ar@{-->}[rr]^{p}&&Q\\
}$$
with all maps and homotopies above is homotopy commutative,
\item if there exists another map $p':P \rightarrow Q$ and homotopies $K':f''\sim p'\circ f'$ and $L': p'\circ g' \sim g''$ such that the obtained diagram is homotopy commutative , then there exists a homotopy $M:p\sim p'$ such that the whole diagram with all maps and homotopies is homotopy commutative.
\end{enumerate}
\end{definition}

There is a ``standard'' construction of the homotopy pushout of any two maps $f: X \rightarrow Y$, $g: X \rightarrow Z$ as: $P_{(f,g)}\cong (Y \amalg  (X \times I)  \amalg  Z) /  \sim$, where $I$ is the unit interval and the equivalence relation $\sim$ is defined as follows. For every $x \in X$, $(x \times 0)$ is identified to $f(x) \in Y$ and $(x \times 1)$ is identified to $g(x) \in Z$.

Up to now, we was working with standard topology, from now on we switch to the simplicial framework. To compute the homology groups of $P_{(f,g)}$, we use the effective homology method, which often allows one to determine the homology groups of an initial space $S$ by building a homotopy equivalence between the associated chain complex $C_\ast(S)$ (it is worth noting that $H_\ast(S) = H_\ast(C_\ast(S))$) and an \emph{effective} (of finite type) chain complex $E_\ast(S)$. This homotopy equivalence is denoted by $C_\ast(S) \lrrdc E_\ast(S)$. The homology groups $H_\ast(E_\ast(S))$ are easily computable and the equivalence provides an isomorphism $H_\ast(S)\cong H_\ast(E_\ast(S))$ which makes it possible to determine the looked-for homology groups of $S$. We urge the interested reader to consult a description of this technique in \cite{RS02}.

In our case, given $f:X\rightarrow Y$ and $g:X\rightarrow Z$ where $X,Y$ and $Z$ are simplicial sets with effective homology, we must build an effective homology version of $P_{(f,g)}$ (that is, an equivalence between $C_\ast(P_{(f,g)})$ and an effective chain complex $E_\ast(P_{(f,g)})$. The construction of the effective homology of $P_{(f,g)}$ needs several steps, involving, for instance, the building of the cone of morphisms and the suspension of chain complexes. A complete description of our construction of the effective homology of $P_{(f,g)}$ can be found in \cite{Her09}. Thus we have designed an algorithm with the following input-output description.

\begin{Alg}\label{alg1}
\inp two simplicial morphisms $f:X\rightarrow Y$ and $g:X\rightarrow Z$ where $X,Y$ and $Z$ are simplicial sets with effective homology.\\
\outp the effective homology version of $P_{(f,g)}$, that is, an equivalence $C_\ast(P_{(f,g)}) \lrrdc E_\ast(P_{(f,g)})$, where $E_\ast(P_{(f,g)})$ is an effective chain complex.
\end{Alg}

\section{The algorithm}

Let $X, Y$ and $Z$ be simplicial sets with effective homology, and, given $f:X\rightarrow Y$ and $g:X\rightarrow Z$ simplicial morphisms, we can define algorithmically the effective homology version of the pushout as follows.

\begin{minipage}[m]{.46\textwidth}
{\scriptsize
Let $f:X\rightarrow Y$ and $g:X\rightarrow Z$:
\begin{enumerate}
  \item construct\\ $rc:= (X\times I)\backslash (X\times \{0\} \cup X\times \{1\})$,
  \item construct $ds:= Y\oplus Z$,
  \item construct $sds:=$ the suspension of $ds$,
  \item construct the pushout $p:=P_{(f,g)}$,
  \item construct the morphisms\\ $rc \xrightarrow{\sigma} p \xrightarrow{\rho} ds \xrightarrow{shift} sds$,
  \item construct $\chi:=shift \circ \rho \circ \sigma$,
  \item construct $cone:=$ cone of $\chi$,
  \item construct the effective homology version of\\ $cone:=$ $Cone(\chi)\lrrdc ECone(\chi)$,
  \item construct the morphisms\\ $\xymatrix{P_{(f,g)}\ar@<0.5ex>[r]^{f}&Cone(\chi)\ar@<0.5ex>[l]^{g}}$,
  \item from $\xymatrix{P_{(f,g)}\ar@<0.5ex>[r]^{f}&Cone(\chi)\ar@<0.5ex>[l]^{g}}\lrrdc ECone(\chi)$\\ construct $P_{(f,g)}\lrrdc ECone(\chi)$.
\end{enumerate}}
\end{minipage}
\begin{minipage}[m]{.005\textwidth}

\end{minipage}
\fbox{\begin{minipage}[m]{.485\textwidth}

\bmpi \textbf{Input:} $f:X\rightarrow Y$ \boxtt{, } $g:X\rightarrow Z$ \empi
\bmpi \texttt{\textcolor{white}{aaa}}$sorc\leftarrow$ \texttt{sorc(}$f$\texttt{)} \empi
\bmpi \texttt{\textcolor{white}{aaa}}$trgtf\leftarrow$ \texttt{trgt(}$f$\texttt{)}\empi
\bmpi \texttt{\textcolor{white}{aaa}}$trgtg\leftarrow$ \texttt{trgt(}$g$\texttt{)}\empi
\bmpi \texttt{\textcolor{white}{aaa}}$rc\leftarrow$ \texttt{remove-covers(}$sorc$\texttt{)}\empi
\bmpi \texttt{\textcolor{white}{aaa}}$ds\leftarrow$ \texttt{direct-sum(}$trgtf$,$trgtg$\texttt{)} \empi
\bmpi \texttt{\textcolor{white}{aaa}}$sds\leftarrow$ \texttt{suspension-functor(}$ds$\texttt{)} \empi
\bmpi \texttt{\textcolor{white}{aaa}}$p\leftarrow$ \texttt{pushout(}$f$,$g$\texttt{)} \empi
\bmpi \texttt{\textcolor{white}{aaa}}$\sigma\leftarrow$ \texttt{build-mrph(}$rc$,$p$,$0$,$rc\rightarrow p$\texttt{)} \empi
\bmpi \texttt{\textcolor{white}{aaa}}$\rho\leftarrow$ \texttt{build-mrph(}$p$,$ds$,$0$,$p\rightarrow ds$\texttt{)} \empi
\bmpi \texttt{\textcolor{white}{aaa}}$shift\leftarrow$ \texttt{build-mrph(}$ds$,$sds$,$1$,$ds\rightarrow sds$\texttt{)} \empi
\bmpi \texttt{\textcolor{white}{aaa}}$\chi\leftarrow$ \texttt{build-mrph(}$rc$,$sds$,$0$,$rc\rightarrow sds$\texttt{)} \empi
\bmpi \texttt{\textcolor{white}{aaa}}$cone\leftarrow$ \texttt{cone2(}$\chi$\texttt{)} \empi
\bmpi \texttt{\textcolor{white}{aaa}}$cone$-$efhm \leftarrow$ \texttt{cone-efhm(}$cone$\texttt{)}\empi
\bmpi \texttt{\textcolor{white}{aaa}}$f\leftarrow$ \texttt{build-mrph(}$cone$,$p$,$0$,$cone\rightarrow p$\texttt{)} \empi
\bmpi \texttt{\textcolor{white}{aaa}}$g\leftarrow$ \texttt{build-mrph(}$p$,$cone$,$0$,$p\rightarrow cone$\texttt{)} \empi
\bmpi \texttt{\textcolor{white}{aaa}}$lf \leftarrow$ \texttt{cmps(}$f$,\texttt{lf(}$cone$-$efhm$\texttt{)}\texttt{)}\empi
\bmpi \texttt{\textcolor{white}{aaa}}$lg \leftarrow$ \texttt{cmps(}\texttt{lg(}$cone$-$efhm$\texttt{)},$g$\texttt{)}\empi
\bmpi \texttt{\textcolor{white}{aaa}}$lh \leftarrow$ \texttt{lh(}$cone$-$efhm$\texttt{)}\empi
\bmpi \texttt{\textcolor{white}{aaa}}$rf \leftarrow$ \texttt{rf(}$cone$-$efhm$\texttt{)}\empi
\bmpi \texttt{\textcolor{white}{aaa}}$rg \leftarrow$ \texttt{rg(}$cone$-$efhm$\texttt{)}\empi
\bmpi \texttt{\textcolor{white}{aaa}}$rh \leftarrow$ \texttt{rh(}$cone$-$efhm$\texttt{)}\empi
\bmpi \texttt{\textcolor{white}{aaa}}$pushout$-$efhm \leftarrow$ \texttt{build-hmeq(}$lf$,$lg$,$lh$,$rf$,$rg$,$rh$\texttt{)}\empi
\bmpi \texttt{\textcolor{white}{aaa}}\textbf{return}\texttt{(}$pushout$-$efhm$\texttt{)}\empi

\end{minipage}}

The arguments of the form $a\rightarrow b$ represent pure lisp functions implementing the mathematical algorithm of a morphism, see \cite{Kenzo} for a complete description of this kind of functions.

In the above algorithm several functions are used. On the one hand, some of them (\boxtt{sorc}, \boxtt{trgt}, \boxtt{lf}, $\ldots$) are Kenzo functions on morphisms and homotopy equivalences which allow one to access to the internal components of these mathematical structures (for instance, \boxtt{sorc} and \boxtt{trgt} return the source and the target simplicial set of a morphism respectively). On the other hand, functions such as \boxtt{remove-covers} or \boxtt{direct-sum} are new algorithms introduced by us to construct the effective homology version of the pushout of simplicial morphisms. Moreover, in spite of not appearing explicitly in the effective homology version of the pushout algorithm, relevant functions have been implemented for intermediate constructions. Let us focus on some important ones.

Let $\xymatrix{0&A_\ast\ar[l]_0\ar@<0.5ex>[r]^{\sigma}&B_\ast\ar@<0.5ex>[l]^j\ar@<0.5ex>[r]^{\rho}&C_\ast\ar@<0.5ex>[l]^i&0\ar[l]}$ be an effective short exact sequence of chain complexes ($A$, $B$ and $C$ are chain complexes and the morphisms satisfy $\rho i = id_{C_\ast}$, $i\rho + \sigma j = id_{B_\ast}$ and $j\sigma = id_{A_\ast}$, $i$ and $j$ are not compatibility in general with the differentials). Then three general algorithms are available: $SES_1: (B_{\ast,EH}, C_{\ast,EH}) \mapsto A_{\ast,EH}$, $SES_2: (A_{\ast,EH}, C_{\ast,EH}) \mapsto B_{\ast,EH}$ and $SES_3: (A_{\ast,EH}, B_{\ast,EH}) \mapsto C_{\ast,EH}$ producing a version with effective homology of one chain complex when versions with effective homology of both others are given. A complete description of these algorithms can be found in \cite{RS06}. The algorithms $SES_1$ and $SES_2$ play a key role in the construction of the effective homology version of the pushout of simplicial sets but they are not included in the Kenzo system, so we implemented both of them. For instance, given a short exact sequence where $B$ and $C$ are chain complexes with effective homology, then $SES_1$ is algorithmically specified as follows.

\bmp
\bmpi \textbf{Input:} $A$\verb|,| $B$\verb|,| $C$\verb|,| $\rho: B\rightarrow C$\verb|,| $\sigma: A\rightarrow B$\verb|,| $i: C\rightarrow B$ and $j: B\rightarrow A$ \empi
\bmpi \verb|   |$cone\leftarrow$ \verb|cone(|$i$\verb|)| \empi
\bmpi \verb|   |$cone$-$efhm\leftarrow$ \verb|efhm(|$cone$\verb|)| \empi
\bmpi \verb|   |$rrdct\leftarrow$ \verb|rrdct(|$cone$-$efhm$\verb|)| \empi
\bmpi \verb|   |$lrdct\leftarrow$ \verb|lrdct(|$cone$-$efhm$\verb|)| \empi
\bmpi \verb|   |$A$-$cone$-$rdct\leftarrow$ \verb|aibjc-rdct(|$A$\verb|,| $i$\verb|,| $\sigma$\verb|,| $B$\verb|,| $j$\verb|,| $\rho$\verb|,| $C$\verb|)|\empi
\bmpi \verb|   |$final$-$lrdct\leftarrow$ \verb|cmps(|$A$-$cone$-$rdct$\verb|,| $lrdct$\verb|)|\empi
\bmpi \verb|   |$ses1$-$hmeq\leftarrow$ \verb|build-hmeq(|$final$-$lrdct$, $rrdct$\verb|)| \empi
\bmpi \verb|   |\textbf{return}\verb|(|$ses1$-$hmeq$\verb|)|\empi
\emp

Let us note, that both \boxtt{cone} (available in Kenzo) and \boxtt{cone2} (implemented by us) algorithms implement the construction of the cone of a morphism taking two different but equivalent definitions, each one of them useful in a concrete context.

As we claimed at the beginning of this paper, many of the usual constructions in Topology are particular cases of pushout constructions. In this way, algorithms to build, for instance, the effective homology version of wedges or joins are obtained as instances from the pushout algorithm explained in this section. As an example, the $8$-th homology group of the join of $K(\mathbb{Z},2)$ and $K(\mathbb{Z},3)$, that is $H_8(K(\mathbb{Z},2)\ast K(\mathbb{Z},3))$, can be computed as follows.

\bmp
\bmpi\verb|> (homology (join (k-z 2) (k-z 3)) 8)|\empim
\bmpi\verb|Homology in dimension 8:|\empi
\bmpi\verb|Component Z/2Z|\empi
\bmpi\verb|Component Z|\empi
\emp

\section{Conclusions and further work}

The algorithm presented in this paper allows one to build the pushout of $f:X\rightarrow Y$ and $g: X \rightarrow Z$, where $X, Y$ and $Z$ are simplicial sets. The implementation has been written in Common Lisp, as a new module (1600 lines) enhancing the Kenzo system. A more challenging task is the implementation of the dual notion of pushout, called \emph{pullback} \cite{Mat76}, that remains as further work.

\bibliographystyle{plain}
\bibliography{ehotposs}

\end{document}